\documentclass[12pt, a4paper,pdftex]{amsart} 
\makeatletter\let\@wraptoccontribs\wraptoccontribs\makeatother 

\title{Projective hypersurfaces in tropical scheme theory I: the Macaulay ideal} 

\author{Alex Fink}
\author{Jeffrey Giansiracusa}
\author{Noah Giansiracusa}
\contrib{with an appendix by Joshua Mundinger}

\date{\today}

\usepackage[T1]{fontenc}

\usepackage[protrusion=true,expansion=true]{microtype}

\usepackage{hyperref}
\usepackage{amssymb}
\usepackage{mathrsfs}  
\usepackage{latexsym}
\usepackage{amsmath}
\usepackage{color}
\usepackage{enumitem}
\usepackage{graphicx}
\usepackage{blkarray}
\usepackage{bbm} 

\usepackage{mathtools}

\usepackage[cmtip,arrow]{xy}
\usepackage{pb-diagram,pb-xy}
\usepackage{tikz}
\tikzset{bdot/.style={circle, thick, draw=black, top color=black!20, bottom color=black}, 
              wdot/.style={circle, thick, draw=black!20, top color=white, bottom color=white}}

\usepackage[margin=1.0in]{geometry}
\numberwithin{equation}{subsection}

\newtheorem{thmA}{Theorem}

\newtheorem{atheorem}{Theorem}[section]  
\newtheorem{theorem}{Theorem}[subsection]  
\newtheorem{question}[theorem]{Question} 
\newtheorem{lemma}[theorem]{Lemma} 
\newtheorem{alemma}[atheorem]{Lemma} 
\newtheorem{proposition}[theorem]{Proposition}
\newtheorem{aproposition}[atheorem]{Proposition}

\newtheorem{recipe}[theorem]{Recipe}
\theoremstyle{remark} 
\newtheorem{remark}[theorem]{Remark}
\newtheorem{aremark}[atheorem]{Remark}
\newtheorem{definition}[theorem]{Definition}
\newtheorem{adefinition}[atheorem]{Definition}
\newtheorem{example}[theorem]{Example}

\newcommand{\D}{\mathcal{D}}

\newcommand{\B}{\mathbb{B}}

\newcommand{\R}{\mathbb{R}}

\newcommand{\PP}{\mathbb{P}}
\newcommand{\T}{\mathbb{T}}

\newcommand{\proj}{\mathrm{Proj}\:}
\newcommand{\trop}{\operatorname{trop}}

\newcommand{\supp}{\operatorname{supp}}

\newcommand{\rank}{\operatorname{rank}}
\newcommand{\corank}{\operatorname{corank}}

\newcommand{\field}{\mathbbm{k}}

\newcommand{\semifield}{\mathbb{S}}
\newcommand{\fieldpolys}{\mathbbm{k}[\mathbf{x}]}
\newcommand{\semifieldpolys}{\mathbb{S}[\mathbf{x}]}
\newcommand{\Bpolys}{\mathbb{B}[\mathbf{x}]}

\newcommand{\initial}{\mathrm{in}}

\usepackage{setspace}
\usepackage{todonotes}

\begin{document}

\begin{abstract}
  A ``tropical ideal'' is an ideal in the idempotent semiring of tropical polynomials that is also,
  degree by degree, a tropical linear space.  We introduce a construction based on transversal matroids that canonically extends any principal ideal to a tropical ideal.  We call this the Macaulay tropical ideal.  It has a universal property: any other extension of the given principal ideal to a tropical ideal with the expected Hilbert function is a weak image of the Macaulay tropical ideal.  For each $n\geq 2$ and $d\geq 1$ our construction yields a non-realizable degree $d$ hypersurface scheme in $\mathbb{P}^n$. Maclagan-Rinc\'on produced a non-realizable line in $\mathbb{P}^n$ for each $n$, and for $(d,n)=(1,2)$ the two constructions agree.  An appendix by Mundinger compares the Macaulay construction with another method for canonically extending ideals to tropical ideals.  
\end{abstract}

\maketitle


\section{Introduction}

\subsection{Background and motivation} 
A foundational question in tropical geometry is: \emph{what are the fundamental geometric objects of the tropical world?} As the subject has developed, various answers to this question have emerged.  Over the past few years the perspective of \emph{tropical ideals} \cite{Maclagan-Rincon, Maclagan-Rincon2, Maclagan-Rincon3}, and their associated \emph{tropical schemes} \cite{GG1, Lorscheid-hyperfields, Lorscheid-unifying}, has appeared as an intriguing proposed answer.  However, working with tropical ideals has turned out to be much more difficult than working with classical ideals.

Let $\T$ denote the tropical semifield.  A tropical ideal in the algebra $\T[x_0, \dots, x_n]$ of tropical polynomials is an ideal that is also a tropical linear space.  Tropical ideals satisfy the ascending chain condition, admit Hilbert polynomials, have balanced polyhedral complexes as their varieties, and satisfy a tropical Nullstellensatz \cite{Maclagan-Rincon2}.  However, finitely generated ideals in $\T[x_0, \dots, x_n]$ generally fail to be tropical ideals, and tropical ideals generally fail to be finitely generated.  Moreover, examples show that no finite degree truncation suffices to uniquely determine a tropical ideal, and there can in fact be uncountably many extensions of a given truncation to a full tropical ideal \cite{Anderson-Rincon-paving}. Producing tropical ideals is significantly more difficult than producing ideals over a field.  

In this paper we study the problem of extending principal ideals to tropical ideals.  Given a tropical polynomial $f$, we look for a tropical ideal that contains $f$ (and hence the principal ideal generated by $f$). This question is only interesting if the tropical ideal is not allowed to be too large.  Over a field, there is of course a unique minimal ideal containing a given polynomial $F$, namely the principal ideal generated by $F$.  In contrast, in the poset of tropical ideals containing a tropical polynomial $f$ (ordered by inclusion), one often finds many minimal elements.  We propose the following definition, which generalizes the definition given in \cite{Silversmith}; see section \ref{sec:tropically-principal-ideals} for discussion.
\begin{definition} 
  Given a tropical polynomial $f$, a tropical ideal $J$ is said to be \emph{tropically principal over $f$} it is minimal in the poset of tropical ideals containing $f$, and $J$ is said to be a \emph{tropically principal ideal} if there exists an $f$ over which it is tropically principal.
\end{definition}
In the algebra of tropical polynomials, principal ideals are usually not tropical ideals, and tropically principal ideals are usually not principal ideals.

One abundant source of tropical ideals is through tropicalization of classical ideals, but just as not all matroids are realizable, not all tropical ideals are realizable. One can produce an extension of a tropical polynomial $f$ to a tropically principal ideal over $f$ by lifting $f$ to a classical polynomial $F$ and then tropicalizing the principal ideal $\langle F \rangle$ generated by $F$ (note that the result can depend on the field that one lifts to).  Thus, we are most interested in constructing non-realizable tropically principal ideals.

\subsection{Results}
This paper is an exploration of the set of tropically principal ideals over a given tropical
polynomial $f$.  We describe a combinatorial recipe to produce a tropically principal ideal $[f]$.
The construction is built on the operation of \emph{stable sum} of tropical linear spaces,
which is the extension to valuated matroids of the matroid union operation \cite{Fink-Rincon}. 

Given a homogeneous polynomial $f$, for each $d \ge \deg(f)$ we form the Macaulay
matrix $\D_{d}(f)$, which is the matrix whose columns and rows correspond to all monomials
of degree $d$ and $d-\deg(f)$, respectively, and whose entries record the coefficients of the
multiples of $f$ by degree $d-\deg(f)$ monomials. We then consider the tropical linear space 
$[f]_d$ (inside the space of degree $d$ polynomials) given by the stable sum of the lines spanned by the rows of $\D_d(f)$. 

\begin{thmA}\label{thm:A}
  The tropical linear spaces $[f]_d$ form a tropically principal ideal $[f]$ over $f$. We call this the \emph{Macaulay tropical ideal} generated by $f$.  It satisfies the following conditions: 
\begin{enumerate}
\item The associated tropical variety is the tropical hypersurface of $f$.
\item The Hilbert function of $[f]$ is that of a hypersurface of degree $\deg(f)$.
\item If $f$ has at least three terms and each term contains a variable not appearing in any other term, then $[f]$ is not realizable over any field. 
\item When $f=x_0 + x_1 + x_2$, the ideal $[f]$ equals the $\PP^2$ case of the
  non-realizable tropical ideal of a line given in \cite[Example 2.4]{Maclagan-Rincon2}.
\end{enumerate}
\end{thmA}

Section \ref{sec:weak order} shows that
the Macaulay ideal $[f]$ satisfies a universal property with respect to the so-called weak order on matroids. 
If the underlying matroids of a tropically principal ideal $I$ containing $f$
have the same ranks as in~$[f]$ (we call these \emph{numerically principal}),
then these matroids are \emph{weak images} of the corresponding matroids of~$[f]$.  This means that the set of non-zero Pl\"ucker coordinates of $I$ is a subset of the non-zero Pl\"ucker coordinates of $[f]$.

In general, the set of all tropically principal ideals over a given $f$ appears to be complicated, but in the cases of the following theorem, proved in Proposition~\ref{prop:quadbinom}, they are singletons.
\begin{thmA}\label{thmA:B}
Let $f$ be a homogeneous tropical monomial or binomial. 
Then the only tropically principal ideal over $f$ is $[f]$,
which is equal to the realizable ideal $\trop(\langle F \rangle )$ for any lift $F$ of $f$ to a non-archimedean field.
\end{thmA}
A subsequent paper will continue the investigation of the structure of the set of tropically principal ideals over~$f$
when $f$ is a linear form or a homogeneous quadratic in two variables.

An appendix by Joshua Mundinger introduces another method for canonically extending a principal ideal in the algebra of tropical polynomials to a tropical ideal, by iteratively and systematically including circuit eliminations in a minimal way.  When applied to a principal ideal, this construction is not guaranteed to yield a tropically principal ideal, but for binomials it does and, consequently, it coincides with the Macaulay tropical ideal construction in that case.

\subsection{Outline of the paper}

Section \ref{sec:preliminaries} sets up our notation and conventions, reviewing the basics the algebra that we will need, tropical linear spaces, tropical ideals, and some required matroid theory.  Section \ref{sec:tropically-principal-ideals} explores our definition of tropically principal ideals and its relation to a similar definition proposed by Silversmith \cite{Silversmith}.  Finally, in Section \ref{sec:Macaulay} we come to the heart of the paper that defines and studies the tropically principal ideal $[f]$.

\subsection*{Acknowledgements}
We thank Nicholas Anderson, Dori Bejleri, Tyler Foster, Diane Maclagan, Andrea Pulita, and Felipe Rinc\'on for helpful conversations. We also thank the American Institute of Mathematics for hosting the workshop at which this project began, and Diane for catalyzing this paper's completion with the support of EPSRC grant EP/X02752X/1.
NG was supported in part by NSF DMS-1802263. AF was supported by EPSRC grant EP/X001229/1, and by the European Union's Horizon 2020 research and innovation programme under the Marie Sk\l odowska-Curie grant agreement No 792432. 
JG was supported by EPSRC grant EP/Y028872/1.

\section{Preliminaries and notation}\label{sec:preliminaries}

\subsection{Idempotent semifields, tropical polynomials, and non-archimedean seminorms}

Throughout this paper $\semifield = (\semifield, +, \cdot)$ denotes a totally ordered idempotent
semifield; this means $\semifield \setminus\{0\}$ is a totally ordered abelian group, written
multiplicatively, and $\semifield$ is equipped with the maximum operation for its addition. In
particular, the additive and multiplicative units are $0$ and $1$, respectively, and $0$ is the
minimal element of $\semifield$.  Tropical geometry typically works with coefficients in the extended real line
$\T=\R\cup \{\infty\}$ with operations $(\min,+)$, but there have been some efforts to develop higher rank tropical geometry
\cite{Dhruv-Tyler-Hahn} and most of our constructions and proofs extend immediately to this more
general setting.  The Boolean semifield $\mathbb{B} = \{0,1\}$ is the unique 2-element idempotent
semifield; it is a subsemifield of $\semifield$.

We write $\semifieldpolys$ for the set of polynomials on the variables $\mathbf{x}=\{x_0 \dots, x_n\}$ with coefficients in $\semifield$.  This has the structure of an idempotent semiring, and the unit morphism $\semifield \to \semifieldpolys$ makes it an $\semifield$-algebra.

To avoid notational confusion, we will work with \emph{non-archimedean seminorms} rather than valuations.
Let $R$ denote a ring.  A non-archimedean seminorm is a map $\nu: R \to \semifield$ 
such that $\nu(0) = 0$, $\nu(1) = \nu(-1) = 1$, and for all $a,b\in R$,
\begin{enumerate}
\item $\nu(ab) = \nu(a)\nu(b)$, and
\item $\nu(a+b) \leq \nu(a) + \nu(b)$. 
\end{enumerate}
A seminorm is a \emph{norm} if $\nu^{-1}(0) = \{0\}$. Note that the distinction between seminorms
and valuations is purely in how we view $\semifield$, namely whether we view $\semifield^\times$ as a
multiplicative or additive group.  The exponential map  $x \mapsto -e^x$ gives an isomorphism 
\[
(\R\cup \{\infty\}, \mathrm{min}, +) \cong (\R_{\geq 0}, \mathrm{max},\times).
\]
By working with seminorms we avoid the usual confusion in tropical algebra where $0$ is
the multiplicative identity and $\infty$ is the additive identity.

Throughout this paper we assume that $\nu: \field \to \semifield$ is an algebraically closed field
equipped with a surjective non-archimedean norm (if $\nu$ is not surjective then one can replace
$\semifield$ with the image of $\nu$ \cite[Lemma 2.5.3]{GG1}, or pass to an extension of $\field$).
Under these hypotheses, one can choose a multiplicative section
\[
s: \semifield \to \field
\]
of the norm. This is proven for $\T$ in \cite[Lemma 2.1.15]{Maclagan-Sturmfels} and the proof
extends immediately to any totally ordered $\semifield$.  

We say that a sum of elements of $\semifield$ \emph{tropically vanishes} if the value is unchanged by dropping any single term, equivalently, the maximum in the sum is attained at least twice.

\subsection{Tropical orthogonal duals}

Given elements $f = \sum_i f_i e_i$ and $g=\sum_i g_i e_i$ in $\semifield^n$, we say that they are \emph{tropically orthogonal} if the sum
\[
\sum_i f_i g_i \in \semifield
\]
tropically vanishes.
Given a subset $K \subset \semifield^n$, the tropical orthogonal dual, denoted $K^\perp$, is the set 
\[
\{v \in \semifield^n \quad | \quad \text{$v$ is tropically orthogonal to $x$ for all $x\in K$} \}.
\]

\subsection{Valuated matroids and tropical linear spaces}\label{sec:TLS}

Valuated matroids were introduced first in \cite{Dress-Wenzel-1, Dress-Wenzel-2}, then later from a
different perspective in \cite{Speyer}.  We recommend \cite[\S4.1]{Frenk} for an excellent
introduction and algebraic treatment.  The connection with the classical Pl\"ucker embedding is
explained in \cite{GG3}. Here we provide a quick review of the basic definitions that will be used
in this paper.

Let $E$ be a finite set, and let $\binom{E}{d}$ denote the set of size $d$ subsets of $E$.  A
nonzero point $p \in \semifield^{\binom{E}{d}}$ is said to be a rank $d$ \emph{valuated matroid} or
\emph{tropical Pl\"ucker vector} if for any $A\in \binom{E}{d+1}$ and $B\in \binom{E}{d-1}$, the sum
\begin{equation}\label{eq:Plucker}
\sum_{i\in A\smallsetminus B} p_{A \smallsetminus \{i\}}\,p_{B\cup \{i\}} \in \semifield
\end{equation}
tropically vanishes.
The entries of $p$ are called \emph{Pl\"ucker coordinates}.

Given a valuated matroid $p\in \semifield^{\binom{E}{d}}$, for each set $C \in \binom{E}{d-1}$ there is an
associated \emph{cocircuit}
\[
\beta_C = \sum_{i\notin C} p_{C\cup \{i\}}\,e_i \in \semifield^E.
\]
Similarly, for each $D\in
\binom{E}{d+1}$ there is an associated linear form
\[
\alpha_D = \sum_{i\in D}p_{D\smallsetminus \{i\}}\,e_i
\]
called the \emph{fundamental circuit} of $D$.  The Pl\"ucker relations say that the circuits and cocircuits are tropically orthogonal.  In fact, it turns out that 
\[
\{\text{circuits}\}^\perp = \mathrm{span} \{\text{cocircuits}\},  \quad \text{and} \quad \{\text{cocircuits}\}^\perp = \mathrm{span} \{\text{circuits}\}
\]
The \emph{tropical linear space} $L_p \subset \semifield^n$ associated with $p$ is the submodule spanned by the cocircuits.  It can equivalently be described as the intersection of the tropical hyperplanes $\alpha_D^\perp$ defined by the circuits of $p$.  Since $L_p$ uniquely determines $p$ up to a scalar (see \cite{Murota-circuits,Speyer,GG3} for three different perspectives), we often work with $p$ and $L_p$ interchangeably.

Tropical linear spaces can also be characterized intrinsically using a valuated version of the
(co)circuit elimination axiom.  We state here the variant applying to sums of cocircuits (sometimes
called covectors) since this is a bit more natural from the module-theoretic perspective: a
submodule $L\subset \semifield^E$ is a tropical linear space if and only if for any pair $v=\sum v_i
e_i, w=\sum w_ie_i \in L$ with $v_j = w_j \ne 0$ for some $j \in E$, there exists $u=\sum u_i e_i
\in L$ with $u_j=0$ and for all $i\ne j$ the inequality $u_i \le v_i+w_i$ holds with equality for
those $i$ such that $v_i \ne w_i$ \cite[Theorem 3.4]{Murota-circuits}.

A homomorphism of semifields $\semifield \to \semifield'$ sends valuated matroids over $\semifield$
to valuated matroids over $\semifield'$.  An ordinary matroid is the same as a valuated matroid over
the booleans $\mathbb{B} = \{0,1\}$.  Any idempotent semifield $\semifield$ admits a unique
homomorphism $\semifield \to \mathbb{B}$ given by sending all nonzero elements to $1$, and this
homomorphism sends a valuated matroid over $\semifield$ to its \emph{underlying matroid}.

\subsection{Stable sum of valuated matroids}

The \emph{stable sum} operation takes two valuated matroids $p$ and $q$ on a given ground set and,
under mild conditions, produces a new valuated matroid $p\wedge q$ with rank equal to
$\mathrm{rank}(p)+\mathrm{rank}(q)$. It is the extension to valuated matroids of the matroid union operation. The dual operation of stable intersection of tropical linear spaces was first defined by Speyer \cite{Speyer}, and stable sum first appeared in \cite{Frenk} and \cite{Fink-Rincon}.  
The stable sum defines a tropical linear space containing both $L_p$ and $L_q$;
we warn the reader there might not be a unique minimal common superspace of this form.
Given tropical Pl\"ucker vectors $p \in \semifield^{\binom{E}{d}}$ and $q \in
\semifield^{\binom{E}{e}}$, the components of $p\wedge q \in \semifield^{\binom{E}{d+e}}$ are defined by
\[
(p\wedge q)_I = \sum_{A\sqcup B = I} p_A q_B,
\]
where the sum is over disjoint sets $A$ and $B$ of size $d$ and $e$, respectively, whose union is
$I$.  If there is at least one nonzero term in this sum then $p\wedge q$ is a valuated matroid \cite[Proposition 5.1.2]{GG3}.

\subsection{Transversal matroids}\label{sec:TransMat}

A \emph{set system} on a set $E$ is a collection of subsets
$\{A_i \subset E\}_{i \in I}$.  A \emph{partial transversal} for this set system is a subset
$Y \subset E$ such that there exists an injective map $\varphi: Y \hookrightarrow I$ satisfying
$y \in A_{\varphi(y)}$ for each $y\in Y$.  The classical theorem of Edmonds and Fulkerson
\cite{Edmonds-Fulkerson} says that the partial transversals of a finite set system constitute the
independent sets of a matroid.  Matroids arising in this way are called \emph{transversal
  matroids}. All transversal matroids are representable over $\mathbb{Q}$, and over
sufficiently large fields of any finite characteristic \cite{Piff-Welsh}.

A matroid of rank $d$ is transversal if and only if it can be written as the stable sum of $d$
matroids of rank 1. Indeed, a rank $d$ transversal matroid can always be presented by a set system of size $d$ \cite[Theorem 2.6]{BoninTransversalNotes}, say $\{A_i \subset E\}_{i \in I}$ with $|I|=d$, and then it is equal to the matroid union (equivalently, stable sum) of the $d$ rank 1 matroids where the bases of the $i^{th}$ matroid are the elements of $A_i$.  Replacing $\mathbb{B}$ with the semifield $\semifield$, there is then an obvious extension of the
definition of being transversal from matroids to valuated matroids: a valuated matroid is said to be
transversal if it can be written as a stable sum of rank 1 valuated matroids.  The structure of
tropical linear spaces associated with transversal valuated matroids was first studied in
\cite{Fink-Rincon}, where these objects were called \emph{Stiefel tropical linear spaces}.

Here is another description of transversal valuated matroids.  For $k<n$, the rows of a matrix
$N\in \mathrm{Mat}_{k\times n}(\semifield)$ each span a line in $\semifield^n$, and the
submodule generated by these is the rowspace of $N$.  This rowspace will often fail to be a tropical
linear space, and there need not be a unique minimal tropical linear space containing it.
However, if $N$ has a nonzero maximal minor then the stable sum of the rows defines a
canonical rank $k$ tropical linear space containing the rowspace; the corresponding valuated matroid is given by the maximal minors of $N$.  (The situation where $N$ is allowed to have all of its maximal minors equal to zero is studied in \cite{Mundinger}.)

\subsection{Monomials}\label{ssec:monomials}

We identify the set of integral lattice points in the positive orthant $\mathbb{R}^{n+1}_{\geq 0}$ with the set of monomials in variables $x_0, \ldots, x_n$. For $d\geq 1$ let $\Delta_d^n$ denote the set of non-negative integral points of the hyperplane $x_0 + x_1 + \cdots + x_n = d$. By convention, we define $\Delta_0^n = \{1\}$.  This set $\Delta^n_d$ is then identified with the set of degree $d$ monomials in $n+1$ variables.  It has the shape of an $n$-simplex with the pure powers $x_i^n$ at each of the vertices.  Along each edge there are $d+1$ points. The number $N_d^n :=|\Delta_d^n|$ of monomials of degree $d$ is the binomial coefficient $\binom{n+d}{d}$.

For $\ell \le d$, let $\Delta^n_{\ell \to d}$ denote the set of all translates of $\Delta^n_{\ell}$ in $\Delta_d^n$.
There is a canonical bijection
\[
\Delta^n_{\ell \to d} \cong \Delta^n_{d-\ell} 
\]
given by sending a set $A\in \Delta^n_{\ell \to d}$ to its greatest common divisor.

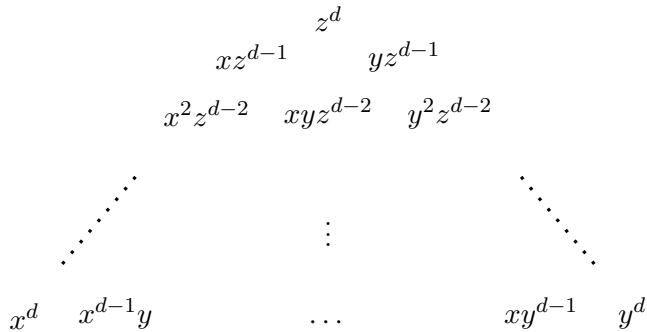
\begin{figure}[hbt]
\begin{center}
\begin{tikzpicture}[font=\small]
  \draw[very thick, loosely dotted] (0.5,1) -- (1.5,2.2); 
  \draw[very thick, loosely dotted] (7.5,1) -- (6.5,2.2); 
  \path (0,0) -- (8,0)
                    node[pos=0,above] {$x^d$}
                    node[pos=0.15,above] {$x^{d-1}y$}
                    node[pos=0.5,above] {$\cdots$}
                    node[pos=0.85,above] {$xy^{d-1}$}
                    node[pos=1,above] {$y^d$};
  \path (0,0) -- (4,5)
                    node[pos=0.6] {$x^2z^{d-2}$}
                    node[pos=0.75] {$xz^{d-1}$};
  \path (8,0) -- (4,5)
                    node[pos=0.6] {$y^2z^{d-2}$}
                    node[pos=0.75] {$yz^{d-1}$};
  \path (4,0) -- (4,5)
                    node[pos=0.3] {$\vdots$}
                    node[pos=0.6] {$xyz^{d-2}$}
                    node[pos=0.85] {$z^d$};
\end{tikzpicture}
\end{center}
\caption{Illustration of the set $\Delta_d^2$.}
\label{fig: Delta_d^2}
\end{figure}

\begin{example}
For $n=2$, the set of monomials 
\[
\{z^d, xz^{d-1},  x^2z^{d-2}, yz^{d-1}, y^2z^{d-2}, xyz^{d-2} \}
\]  
divisible by $z^{d-2}$ is the upper triangle $z^{d-2}\cdot\Delta^2_2$ in $\Delta_d^2$
in Figure~\ref{fig: Delta_d^2}.
\end{example}

\begin{lemma}\label{lem:monomial multiples rank bound}
If a tropical linear space $L\subseteq\semifield^{\Delta^n_d}$ contains $r$ distinct monomial multiples of some polynomial $f\in\semifieldpolys$,
then $\rank L\ge r$.
\end{lemma}

\begin{proof}
Choose any monomial order on $\semifieldpolys$. 
Let $m_1f,\ldots,m_rf$ be the multiples of~$f$ in~$L$,
indexed in increasing order of their greatest term.
The support of any element of~$L$ is a cocyclic set of its underlying matroid $M$.
Cocyclic sets are closed under unions, so $\bigcup_{i_1}^k\supp(m_if)$ is a strictly increasing chain of cocyclic sets of the underlying matroid of~$M$ as $k$ runs from 1 to~$r$. 
Complements of cocyclic sets are flats, so $M$ has a chain of~$r$ flats, implying its rank is at least~$r$.
\end{proof}

\subsection{Tropical ideals}

The degree $d$ piece of $\semifieldpolys$ is the free $\semifield$-module generated by the set
$\Delta^n_d$ of degree $d$ monomials.  It thus makes sense to ask whether a submodule
$L\subset \semifieldpolys_d \cong \semifield^{N^n_d}$ is a tropical linear space.  

\begin{definition}
A homogeneous ideal $J \subset \semifieldpolys$ is a \emph{tropical ideal} if each graded
piece $J_d \subset \semifieldpolys_d$ is a tropical linear space.
\end{definition}

\begin{remark}\label{rem:conventions}
An important note on conventions.  In the convention of \cite{Maclagan-Rincon2}, the degree $d$ part of a tropical ideal $J$ is defined by a valuated matroid $M_d$ such that the elements of $J_d$ are the $\semifield$-linear span of the circuits of $M_d$; thus the rank of $M_d$ is equal to the codimension of $J_d$.

The Hilbert function of a tropical ideal $J$ is the function 
\[
HF_J(d) = \operatorname{corank}(J_d).
\]
By \cite[Corollary 3.8]{Maclagan-Rincon2}, the Hilbert function agrees with a polynomial for all $d$ sufficiently large.

Throughout this paper, we find it more convenient to use the opposite convention: the valuated matroid we associate with $J_d$ has rank equal to the dimension of $J_d$, and $J_d$ is spanned by the cocircuits of this valuated matroid.  This convention is more compatible with the perspective of this paper because we want to think about building $J$ up from elements inside it.
\end{remark}

A few words about the history and context of tropical ideals helps to motivate our investigation of them in this paper.  The scheme-theoretic tropicalization developed in \cite{GG1} proceeds as follows:  a closed subscheme of projective space is defined by a homogeneous ideal $I \subset \fieldpolys$, which can be viewed as a sequence of linear subspaces $I_d \subset \fieldpolys_d$, and by \cite[Proposition 6.1.1]{GG1} the tropical linear spaces $\trop(I_d)$ form an ideal in $\semifieldpolys$.  The tropicalization of the projective subscheme $\proj \fieldpolys /I$ is defined as Proj of a quotient of $\semifieldpolys$ constructed from the tropical ideal $\trop(I)$; but Maclagan and Rinc\'on subsequently proved in \cite{Maclagan-Rincon} that one can reconstruct the tropical ideal $\trop(I)$ from this quotient, and so rather than working with the less familiar congruences used to define quotients in semiring theory, one can stick with a more familiar object, an ideal.  The key property underlying this result of Maclagan and Rinc\'on is not that the ideal is a tropicalization, but that the ideal is a tropical ideal---in other words, to reconstruct the ideal from the quotient it defines relies on properties of tropical linear spaces that do not hold for more general submodules.

The next development, also provided by Maclagan and Rinc\'on, was to show that tropical ideals are not just a convenient way of working with a certain class of projective subschemes that includes the tropicalization of any projective subscheme over a field---they appear to be the ``right'' object to work with since they and the projective subschemes they define have many properties familiar from algebraic geometry that simply do not hold for arbitrary closed subschemes of $\PP^n_\semifield$ \cite{Maclagan-Rincon2}.  In particular, they satisfy the ascending chain condition and admit Hilbert polynomials.

On the other hand, the world of tropical ideals differs significantly  from that of classical ideals.  While classical ideals in $\fieldpolys$ are finitely generated, tropical ideals are generally not finitely generated.  For example, any tropical ideal containing $x+y$ is not finitely generated because it contains the principal ideal $\langle x+y\rangle$ together with a new cocircuit elimination $x^d + y^d$ appearing in each degree $d$.

Maclagan and Rinc\'on showed \cite[Example 3.11]{Maclagan-Rincon2} that 
there is no bound $D$ depending only on the Hilbert function of a tropical ideal $I$ for which the homogeneous parts $(I_i)_{i\leq D}$ of degree at most $D$ determine the whole tropical ideal $I$.  Still more unsettling, Anderson and Rinc\'on showed \cite[Proposition 4.2]{Anderson-Rincon-paving} that there are uncountably many homogeneous tropical ideals in $\mathbb{B}[x,y]$ with constant Hilbert function 3.

\section{Tropically principal ideals}\label{sec:tropically-principal-ideals}

Given a tropical polynomial $f\in \semifieldpolys$, since the principal ideal $\langle f \rangle$ is generally
\emph{not} a tropical ideal if $f$ is not a monomial, we must re-think what it means for a tropical ideal to be ``principal.''
In this section we discuss two alternative definitions of principal for tropical ideals.

\subsection{Two definitions of principal}

Silversmith \cite{Silversmith} defined \emph{tropically principal ideals} in terms of the Hilbert
function.    Here we propose a different definition in terms of the poset of tropical ideals.  We
will show that tropically principal ideals in Silversmith's sense are also tropically principal in
our sense, and we conjecture that the two definitions are in fact equivalent.  For clarity, in this
paper we will work with the next definitions.

\begin{definition}
A \emph{numerically principal tropical ideal} is a homogeneous tropical ideal in $\semifieldpolys$ that has the Hilbert function of a degree $d$ hypersurface for some $d$ (this is the definition that Silversmith gave).  

Given a tropical polynomial $f$, a \emph{tropically principal ideal over $f$} is a minimal element in the poset of tropical ideals containing $f$ under ideal containment, and a \emph{tropically principal ideal} is a tropical ideal that is tropically principal over some $f$.
\end{definition}

We believe it's of interest to study the poset of tropical ideals under containment more generally.
Though we will not use such a generalization in this paper, our definition of ``tropically principal'' generalizes readily to sets of multiple tropical polynomials,
while it's unclear how one should specify the minimal numerics in a generalization of ``numerically principal'' to that setting.

\begin{lemma}\label{lem:tropicallyprincipal}
If $I \subset \semifieldpolys$ is a numerically principal tropical ideal then it is also tropically
principal. More precisely, if $I$ contains $f\in\semifieldpolys$ and has the Hilbert function of a
classical hypersurface of degree $\deg(f)$, then it is tropically principal over $f$.
\end{lemma}

\begin{proof}
Let $I$ be a tropical ideal containing $f$, let $I'$ be a tropical ideal satisfying $f \in I'
\subset I$, and let $F\in\fieldpolys$ be a monomial with $\deg(F) = \deg(f)$.  Our hypothesis is
that $HF_I = HF_{\langle F\rangle}$, and we need to show that $I'=I$.

The containment $I' \subset I$ implies the inequality $HF_{I'} \ge HF_{I}$ of Hilbert functions in
all degrees, and if this inequality is an equality in all degrees then we must have $I=I'$, since in
general a containment of tropical linear spaces is an equality if the ranks coincide.  By applying
the reduction semiring homomorphism $\semifieldpolys \twoheadrightarrow \Bpolys$ we can assume
$\semifield = \B$ and thus can apply \cite[Corollary 3.6]{Maclagan-Rincon2} which says that $HF_{I'}
= HF_{\initial_w(I')}$ for any weight $w\in \mathbb{R}^{n+1}$.  By \cite[Lemma
3.7]{Maclagan-Rincon2} we can choose $w$ so that $\initial_w(I')$ is generated by monomials, in
which case $HF_{\initial_w(I')} \le HF_{\langle F \rangle}$.  Then
\[HF_{I} \le HF_{I'} = HF_{\initial_w(I')} \le HF_{\langle F \rangle} = HF_{I},\]
which establishes the desired equality.  
\end{proof}


\begin{proposition}
If $F\in\fieldpolys$ is homogeneous, then the tropicalized principal ideal $\trop(\langle F \rangle)$ is a numerically principal, and hence tropically principal, over the coefficient-wise valuation $\trop(F)$. 
\end{proposition}

\begin{proof}
This follows from Lemma \ref{lem:tropicallyprincipal} and the observation in \cite[\S7.1]{GG1} that an ideal and its tropicalization have the same Hilbert function.
\end{proof}

We leave the converse of Lemma~\ref{lem:tropicallyprincipal} open.

\begin{question}\label{q:tropical implies numerical}
If $I\subset\semifieldpolys$ is tropically principal,
is it necessarily numerically principal?
\end{question}

\subsection{Tropically principal over a monomial or binomial}

\begin{proposition}\label{prop:quadbinom}
If $f$ is a monomial or a homogeneous binomial, then there is a unique tropically principal ideal over $f$, and it is also numerically principal.  In the monomial case, this ideal is given simply by $\langle f \rangle$.  In the binomial case, the degree $d$ part of this tropically principal ideal is given by the double tropical orthogonal dual $\langle f\rangle_d^{\perp \perp}$.
\end{proposition}

\begin{proof}
If $f$ is a monomial then we can choose a lift to a field, $F\in\fieldpolys$, and we have $\langle f\rangle = \trop(\langle F\rangle)$, from which the conclusion follows.

Suppose $f= f_1 + f_2$ is a binomial.  We can rescale the variables so that both the coefficients are 1, and since this is a monomial change of coordinates it does not affect the number of tropical ideals containing $\langle f\rangle$. Let $I=\langle f\rangle$ be the semiring ideal generated by~$f$.  The set $I_d$ is spanned by the monomial multiples $hf$ for $h \in \Delta^n_{d-\deg f}$. So an element $g=\sum_{m \in\Delta^n_d} g_m m$ lies in $I_d^\perp$ if and only if $g_{hf_1}=g_{hf_2}$ for every such monomial $h$. Let $P_d$ be the finest partition of $\Delta^n_d$ such that $hf_1$ and $hf_2$ lie in the same part of~$P_d$ for every~$h \in \Delta^n_{d-\deg f}$. Then $I_d^\perp$ is spanned as an $\semifield$-module by $\sum_{m\in U}m$ for all parts $U$ of~$P$. This is the tropical linear space associated to the partition matroid of~$P$ (with capacity 1 on all parts) viewed as a valuated matroid with trivial valuation. Finally, $g$ lies in $(I_d)^{\perp\perp}$ if and only if it satisfies the bend relations
\[
\sum_{m\in U}g_m = \sum_{m\in U\setminus u}g_m
\]
for each $U\in P_d$ and $u\in U$. For any $i=0,\ldots,n$ and $U\in P_d$, $x_iU$ is a subset of a single part of~$P_{d+1}$, so $g\in (I_d)^{\perp\perp}$ implies $x_ig\in (I_{d+1})^{\perp\perp}$.

To finish, let $I'$ be any tropical ideal containing~$I$. For any set $S\subseteq\semifield^N$ we have $S\subseteq S^{\perp\perp}$ (because orthogonal dual is a Galois connection), and orthogonal dual is involutive on tropical linear spaces, so
\[
I_d \subset (I_d)^{\perp\perp} \subset (I'_d)^{\perp \perp} = I'_d.
\]
Hence the tropical ideal $I^{\perp\perp}$ is minimal among those containing~$I$.
\end{proof}

On the other hand, below in Theorem~\ref{thm:nonrealize} we show that, for a broad class of $f$ with at least 3 terms, there are multiple distinct tropically principal ideals over~$f$:  
the Macaulay tropical ideal $[f]$ that is the subject of the next section
is distinct from any realizable ideal $\trop\langle F\rangle$ where $\trop(F)=f$.

\section{The Macaulay tropically principal ideal}\label{sec:Macaulay}

In this section we introduce our construction of a canonical tropically principal ideal $[f]$ over a given homogeneous
polynomial $f \in \semifieldpolys$, establish some basic properties, give a sufficient condition for it to be non-realizable, and relate it to the non-realizable tropical ideal examples constructed in \cite{Maclagan-Rincon2}.

\subsection{Construction of the Macaulay tropical ideal \texorpdfstring{$[f]$}{[f]}}

Given a nonzero homogeneous polynomial $f\in \semifieldpolys$, for each $d \ge \deg(f)$ we construct
a matrix $\D_d(f)$ of size $N^n_{d-\deg(f)}\times N^n_d$ whose rows are indexed by the monomials
$\Delta_{d-\deg(f)}^n$, and whose columns are indexed by
the monomials $\Delta_d^n$.  Given monomials $X$ and $Y$ of degrees $(d-\deg(f))$ and $d$,
respectively, the $(X,Y)$ entry of $\D_d(f)$ is the coefficient of $Y$ in $Xf$.  This matrix
$\D_d(f)$ is usually called a \emph{Macaulay matrix}.  

In order to explicitly write down a Macaulay matrix we need to order the monomials in $\semifieldpolys$; unless otherwise stated we shall use the lexicographic order.  In what follows we will also at times refer to the Macaulay matrices $\D_d(F)$ of a homogeneous polynomial $F\in\fieldpolys$, analogously defined.

\begin{example}\label{ex:Macaulay}
For $f= ax_0 + bx_1$ we have
\[
\D_1(f) = \begin{pmatrix} a & b \end{pmatrix},~
\D_2(f) = \begin{pmatrix} a & b & 0 \\ 0 & a & b\end{pmatrix},~
\D_3(f) = \begin{pmatrix} a & b & 0 & 0 \\ 0 & a & b& 0 \\ 0 & 0 & a & b \end{pmatrix}.
\]
\end{example}

\begin{lemma}\label{lem:Macaulay}
Fix a homogeneous polynomial $f\in\semifieldpolys$ and an integer $d \ge \deg(f)$.
\begin{enumerate}
\item The degree $d$ part of the principal ideal $\langle f \rangle$ equals the row space of the Macaulay matrix $\D_d(f)$.
\item There is a nonzero maximal minor of $\D_d(f)$.
\end{enumerate}
\end{lemma}

\begin{proof}
  The first assertion is obvious, so consider the second.  Due to idempotency of addition in $\semifield$, a maximal minor of any matrix over $\semifield$ is nonzero if and only if the corresponding square submatrix has at least one nonzero entry in each row and column.  The monomial multiples of any fixed monomial term in $f$ yield the such a family of nonzero entries. 
\end{proof}

It follows from this lemma that we can form the stable sum of the rows of the matrix $\D_d(f)$, and that this
stable sum contains $\langle f \rangle_d$.

\begin{definition}
  For $d \ge \deg(f)$ let $[f]_d$ denote the Stiefel tropical linear space defined by the Macaulay matrix
  $\D_d(f)$ (i.e., the stable sum of the rows), and for $d < \deg(f)$ let $[f]_d := 0$.  Denote by
  \[
    [f] := \bigoplus_{d \ge 0} [f]_d \subset \semifieldpolys
  \]
  the graded $\semifield$-submodule defined by these pieces.  We refer to $[f]$ as the \emph{Macaulay tropical ideal}
  generated by $f$.
\end{definition}

This name is justified by Proposition \ref{prop:istropicalideal} below.

\begin{example}
  Continuing with Example \ref{ex:Macaulay}, let $f=ax_0+bx_1$.  Then
  \[
    \langle f \rangle_2 = \mathrm{span} \{ x_0 f, x_1 f\}.
  \]
  Assuming $a,b\ne 0$, the principal ideal $\langle f \rangle$ contains the polynomials $ax_0f = a^2x_0^2 + abx_0x_1$ and
  $bx_1 f = abx_0x_1 + b^2x_1^2$, but it does not contain the polynomial $a^2x_0^2 + b^2x_1^2$, which would be required
  by the cocircuit elimination axiom applied to these two polynomials at the term $abx_0x_1$.  On the other hand, the
  tropical Pl\"ucker vector defined by the matrix $\D_2(f)$ is $(a^2, ab, b^2)$, so $[f]_2$ is spanned by the
  following cocircuits:
  \begin{align*}
    \beta_{x_0^2} & = a^2x_0x_1 + abx_1^2,\\
    \beta_{x_0x_1} & = a^2x_0^2 + b^2x_1^2,\\
    \beta_{x_1^2} & = abx_0^2 + b^2x_0x_1.
  \end{align*}
  The $\semifield$-module spanned by the cubic polynomials obtained by multiplying these three cocircuits by $x_0$ and
  $x_1$ also fails to be a tropical linear space and $[f]_3$ provides the missing cocircuit eliminations.  This behavior
  continues in all degrees, so even in this simple case $[f]$ is not finitely generated as an ideal.
\end{example}

\begin{proposition}\label{prop:transversal}
For $d \ge \deg(f)$, the underlying matroid of $[f]_d$ is the transversal matroid represented by the set system
consisting of all translates of $\supp(f)\subset \Delta^n_{\deg(f)}$ in $\Delta^n_d$.
\end{proposition}

\begin{proof}
This follows immediately from the definitions and the discussion in \S\ref{sec:TransMat}.
\end{proof}

In particular, this says that the independent sets of the underlying matroid of $[f]_d$ are the
subsets of $K \subset \Delta_d^n$  such that there exists an injective map
  \[
  \varphi: K \hookrightarrow \Delta_{d-\deg(f)}^n
  \]
with $m \in \supp (\varphi(m)f)$ for each $m\in K$.  The set $K$ determines a set of columns of the Macaulay matrix $\D_d(f)$, and the image $\varphi(K)$ is a set of nonzero entries with one in each column and no two in the same row.

Here is a slightly more geometric description of the bases in the special case where $f$ has full
support.  A subset $K \subset \Delta^n_{d}$ is a basis if there exists a bijection
$\varphi: K \to \Delta^n_{\deg(f) \to d}$ with $m \in \varphi(m)$ for each $m\in K$.  Note that, more
generally, $K$ is an independent set if there exists such a map $\varphi$ that is injective, and it
is a spanning set if there exists such a map that is surjective.  In particular, when $n=2$ and
$\deg(f)=1$, $K$ is a basis if each element in $K$ can be assigned to a unique triangle
$T \cong \Delta^2_1$ containing it.

\subsection{Basic properties of \texorpdfstring{$[f]$}{[f]}}

\begin{proposition}\label{prop:istropicalideal}
The submodule $[f] \subset \semifieldpolys$ is a tropical ideal.
\end{proposition}
\begin{proof}
  To see that $[f]$ is a tropical ideal, it suffices to show that for any variable $x_i$, multiplication by $x_i$ sends
  $[f]_d$ into $[f]_{d+1}$.  Multiplication by $x_i$ gives a coordinate linear inclusion
  $\semifieldpolys_d \hookrightarrow \semifieldpolys_{d+1}$ that sends the homogeneous polynomial corresponding to each
  row of $\D_d(f)$ to a homogeneous polynomial corresponding to a row of $\D_{d+1}(f)$.  The result
  then follows from the observation that stable sum commutes with coordinate linear inclusions.
\end{proof}

\begin{proposition}\label{prop:HilbSeries}
The Hilbert function of $[f]$ coincides with that of a classical hypersurface of degree $\deg(f)$ in $\mathbb{P}_\field^n$.  Consequently, $[f]$ is numerically principal, and hence it is tropically principal over $f$.
\end{proposition}
\begin{proof}
  Let $F\in\fieldpolys$ be a homogeneous polynomial with $\deg(F) = \deg(f)$.  The Macaulay matrix
  $\D_d(F)$ has full rank, which by Lemma \ref{lem:Macaulay}(2) is the rank of the Stiefel tropical linear
  space associated to $\D_d(f)$, so the dimension of the linear space $\langle F\rangle_d$ is equal to the rank of the tropical linear space $[f]_d$ for all $d$.  This proves the first assertion, and the second assertion then follows from Lemma \ref{lem:tropicallyprincipal}.
\end{proof}

\subsection{Non-realizability of \texorpdfstring{$[f]$}{[f]}}

\begin{theorem}\label{thm:nonrealize}
Let $f \in \semifieldpolys$ be a homogeneous polynomial such that
there is a monomial ideal $M$ generated by a subset of the variables~$\mathbf x$ with $|\supp(f)\setminus M|=3$,
and, letting $\supp(f)\setminus M=\{f_1,f_2,f_3\}$, 
we have that $f_j$ does not divide $f_k^if_\ell^{3-i}$ for any $\{j,k,\ell\}=\{1,2,3\}$ and $0\le i\le 3$.
Then $[f]$ is not realizable.  More precisely, if $J$ is
any realizable tropically principal ideal over $f$, then $J$ and $[f]$ are different in degree
$3\deg(f)$.
\end{theorem}

For example, the conditions on $f$ are met if it has at least three terms 
and each term contains a variable not appearing in any other term.

\begin{proof}
  Let us write
  \[
    f = \sum_{i=1}^r c_i f_i
  \]
  with $c_i\in\semifield^\times$ for the decomposition into an $\semifield$-linear combination of $r$ monomials
  $f_i \in \Delta^n_{\deg(f)}$, 
  with the monomials $f_1, f_2, f_3$ being as specified in the theorem hypotheses,
  and set $d := \deg(f)$.  Suppose there exists an ideal $I \subset \fieldpolys$ that
  tropicalizes to $[f]$.  The fact that the Hilbert function of an ideal is preserved under tropicalization, together
  with Proposition \ref{prop:HilbSeries}, implies that $I$ must be a principal ideal, say $I = \langle F\rangle$.  Now
  $[f]_d$ is the one-dimensional tropical linear space spanned by $f$, so the coefficient-wise seminorm $\trop(F)$ must
  be a constant multiple of $f$.  By rescaling $f$ if necessary we can assume $\trop(F) = f$; indeed, rescaling $f$ has
  the effect of rescaling all the entries of the Macaulay matrices defining $[f]$ by the same constant, which means the
  Pl\"ucker coordinates in each degree get rescaled by a power of this constant, but this does not alter the tropical
  linear spaces they define.  Thus we can write
  \[
    F= \sum_{i=1}^r C_i f_i,
  \]
  where $\nu(C_i)=c_i$, and in particular $C_i \ne 0$.  We will show that the tropicalization of $\langle F\rangle$ does
  not equal $[f]$ by showing that there is a Pl\"ucker coordinate that is zero in $\langle F\rangle_{3d}$ and nonzero in
  $[f]_{3d}$.

  First, we partition the monomials of degree $2d$ as follows:
  \[
    U_1 := \{f_1^2, f_1f_2, f_1f_3, f_2^2, f_2f_3, f_3^2\}, ~U_2 = \Delta^n_{2d} \setminus U_1.
  \]
  And we partition the set of degree $3d$ monomials as follows:
  \[
    V_1:=\{f_1^2f_2, f_1^2f_3, f_1f_2^2, f_2^2f_3, f_1f_3^2, f_2f_3^2\}, ~V_2 := \Delta^n_{3d}\setminus V_1.
  \]
  By the assumptions on divisibility of products of $f_1$, $f_2$ and $f_3$, 
  we see that $|U_1|=|V_1|=6$.  We claim that the Macaulay matrix $\D_{3d}(F)$ has the block
  structure
  \[
    \begin{blockarray}{ccc}
      V_1 & V_2 & \\
      \begin{block}{(cc)c}
        P & Q & U_1 \\
        0 & R & U_2 \\
      \end{block}
    \end{blockarray}\ ,
  \]
  meaning that $mf_i \notin V_1$ for each $i=1,\ldots,r$ and each $m\in U_2$.  This also follows from our assumptions:
  if $mf_i = f_jf_kf_\ell$ for $j,k,\ell\in \{1,2,3\}$ then we cannot have $i>3$,
  as this would imply that the left side but not the right side lies in the prime ideal $M$. 
  Hence, $i\in\{1,2,3\}$, implying without loss of generality that $i=j$, but then $m=f_kf_\ell \in U_1$.

  Using this block structure, we see that any maximal minor of $\D_{3d}(F)$ involving all of the $V_1$ columns
  must be of the form
  \[
    \det(P) \cdot \text{(a maximal minor of $R$)}.
  \]
  Moreover, we can explicitly compute the upper-left block,
  again using our assumption on distinctness of products to produce the zeroes:
  \[
    P=\ \begin{blockarray}{ccccccc}
    f_1^2f_2 & f_1^2f_3 & f_1f_2^2 & f_2^2f_3 & f_1f_3^2 & f_2f_3^2 & \\
    \begin{block}{(cccccc)c}
    C_2 & C_3 & 0 & 0 & 0 & 0 & f_1^2  \\
    C_1 & 0 & C_2 & 0 & 0 & 0 & f_1f_2 \\
    0 & C_1 & 0 & 0 & C_3 & 0 & f_1f_3 \\
    0 & 0 & C_1 & C_3 & 0 & 0 & f_2^2  \\
    0 & 0 & 0 & C_2 & 0 & C_3 & f_2f_3 \\
    0 & 0 & 0 & 0 & C_1 & C_2 & f_3^2  \\
    \end{block}
\end{blockarray}\ .
  \]
  Expanding the determinant of this matrix yields two nonzero terms that are each $C_1^2 C_2^2 C_3^2$ but with
  opposite signs, and so $\det(P) = 0$. Hence any maximal minor that includes all of the $V_1$ columns is zero.  On the
  other hand, since $f$ is obtained from $F$ by replacing each $C_i$ with $c_i$, the Macaulay matrix
  $\D_{3d}(f)$ is obtained from $\D_{3d}(F)$ also by replacing each $C_i$ with $c_i$.  In particular,
  $\D_{3d}(f)$ has the same block structure and any maximal minor involving all of the $V_1$ columns is the
  determinant of the upper-left block $\nu(P)$ times a maximal minor of the lower-right block $\nu(R)$---but now
  $\det(\nu(P)) \ne 0$ since idempotency of addition in $\semifield$ means the two terms $c_1^2c_2^2c_3^3$ when
  expanding out $\det(\nu(P))$ do not cancel each other out.

We therefore will have completed the proof once we show that there is a nonzero maximal minor of the tropical matrix $\nu(R)$.  We noted above that multiplication by $f_1$ yields a map $U_2 \to V_2$, and this is certainly injective.  For
each $m\in U_2$, the entry of $\D_{3d}(f)$ at $(m,mf_1)$ is $c_1$; this exhibits a nonzero term in the expansion of this maximal minor of $\nu(R)$, which by idempotency implies that the maximal minor itself is nonzero.  Thus any Pl\"ucker coordinate of $\langle F\rangle_{3d}$ indexed by a set containing $V_1$ is zero, whereas the tropical Pl\"ucker coordinate of $[f]_{3d}$ indexed by $V_1 \cup (f_1\cdot U_2)$ is nonzero.
\end{proof}

In the preceding proof we explicitly established non-realizability by comparing carefully chosen (tropical) Pl\"ucker coordinates.  In the case that $f$ is linear we are able to provide a more conceptual proof.  Indeed, by construction the matroid underlying $[f]_d$, for any homogeneous polynomial $f\in\semifieldpolys$ and any degree $d$, is transversal, whereas the following result shows that for linear $F\in \fieldpolys$ with at least three terms the matroid associated to $\langle F\rangle_d$, which is the matroid in degree $d$ underlying any tropicalization of the ideal $\langle F\rangle$, is non-transversal for $d \ge 3$.

\begin{theorem}
If $L \in \fieldpolys$ is a linear form with at least $3$ terms then for each $d \ge 3$ the matroid associated to $\langle L\rangle_d$ is non-transversal.
\end{theorem}

\begin{proof}
The matroid associated to $\langle L\rangle_d$ is the data of which maximal minors of the Macauley matrix $\D_d(L)$ are nonzero.  By permuting the variables we can assume without loss of generality that $L=c_0x_0 + \cdots + c_sx_s$ with $c_i \ne 0$; moreover, we can further assume $c_i = 1$ for all $i=1,\ldots, s$ since rescaling these variables corresponds to simultaneously rescaling the rows and columns of $\D_d(L)$, and rescaling the rows of a matrix does not alter the row space while rescaling the columns does not alter the matroid associated to the row space.

We use the following characterization of transversal matroids by Bonin
\cite[Theorem 3.6]{BoninTransversalNotes}.
Denote the set of cyclic flats of a matroid $M$ by $\mathcal{Z}(M)$.
Define an integer function $\beta$ on the set of all flats of~$M$ such that 
\begin{equation}\label{eq:Bonin}
\sum_{\substack{G\in\mathcal Z(M)\\ G\supset F}} \beta(G) = \operatorname{corank}(F).
\end{equation}
There is a unique such function, as these equations are triangular with respect to a linear extension of the containment order of flats.  Explicitly, $\beta$ has the following recursive expression:
\[\beta(F) = r(M) - r(F) - \sum_{\substack{G\in\mathcal Z(M)\\ G \supsetneq F}} \beta(G)\]
Then $M$ is transversal if and only if $\beta(F)\ge0$ for all flats~$F$.

Recall that the flats of a matroid are intersections of hyperplanes, and the hyperplanes are the complements of cocircuits.  Hence the flats of the underlying matroid of $\langle L\rangle_d$ are the complements of supports of elements in ~$\langle L\rangle_d$.  Also, a flat $F\subset\Delta_d^n$ is cyclic if none of the sets $F\setminus\{u\}$ is also a flat for $u\in F$. So for each monomial multiple $mL$ of~$L$ with $\deg mL=d$, the complement $H_m$ of its support is a hyperplane because $mL$ is of minimal support among elements of~$\langle L\rangle_d$. Moreover, $H_m$ is cyclic:  no element of~$\langle L\rangle_d$ has support of the form $\supp(mL)\cup\{m'\}$ for $m'\not\in\supp(mL)$,  because the Newton polytope of $L$ would not be a summand of the resulting Newton polytope. (This is where $s\ge 2$ is used.) Therefore $\beta(H_m)=1$ for each $m\in\Delta_{d - \deg L}^n$, the only other summand in~\eqref{eq:Bonin} being $\beta(\Delta_d^n)=0$. Observe that the number of these $H_m$ is $\rank {\langle L\rangle_d}$.

For each monomial $m$ of degree $d-3$, the tropical linear space
$\langle L\rangle_d$ also contains the polynomial
\[
\left( \sum_{i=0}^s x_i^2 - \sum_{0 \le i<j \le s} x_i x_j \right) \cdot Lm = \left( \sum_{i=0}^s x_i^3 -
  3\sum_{0 \le i<j<k \leq s} x_i x_j x_k \right)\cdot m
\]
whose support is
\[
S_m := \{x_i^3m \}_{i=0}^s \cup \{x_ix_jx_km\}_{0 \le i<j<k \leq s}.
\] 
We claim that the complement of $S_m$ is a cyclic flat,
that is, there is no element of~$\langle L\rangle_d$ whose support is $S_m$ together with one further monomial $m'$.
It is impossible that $m'$ lie outside the convex hull of~$S_m$ because of the Newton polytope, as before.
So $m$ divides~$m'$, and we may write a putative element of $\langle L\rangle_d$ with this support
as $mfL$ where $\deg f=2$.
Now $m' = x_i^2x_jm$ for some distinct $i,j\in\{0,\ldots,s\}$.
The condition that $x_{i'}^2x_{j'}m$ is not in the support of~$mfL$ for any pair of distinct indices $(i',j')\neq(i,j)$
implies that the coefficient in $f$ of~$x_{i'}^2$ is the negative of the coefficient of $x_{i'}x_{j'}$.
But all of these equations imply the same equation for $(i',j')=(i,j)$, 
which contradicts the presence of $m'$ in the support of~$mfL$.

So we have found further cyclic flats $F_m := \Delta_d^n\setminus S_m$.  Consider one such.  It is not contained in any
$H_{m'}$, so in order for \eqref{eq:Bonin} to hold at $F_m$ there must be a cyclic flat $G\supset F_m$ with
$\beta(G)\ge1$.  Now the positive summands on the left of~\eqref{eq:Bonin} for the flat $\emptyset$ sum to a value
greater than $\corank\emptyset = \rank \langle L\rangle_d$, so there must be a negative summand as well.  This completes
the proof.
\end{proof}

\begin{remark}
While the preceding theorem gives an alternative, and more conceptual, proof of non-realizability of the tropical ideal $[f]$ in the case of a linear form, it is noteworthy that this theorem statement itself is ``classical'' in the sense that it does not involve any tropical geometry. 
It is a statement about the matroidal structure of a principal ideal and thus concerns a topic that could have been explored decades ago but appears, to our knowledge, not to have previously appeared in the literature.
\end{remark}

\subsection{The tropical variety associated with \texorpdfstring{$[f]$}{[f]}}

Given a homogeneous ideal $J \subset \semifieldpolys$, there are two equivalent ways of defining the associated tropical
variety $V^{\trop}(J) \subset \PP^n(\semifield)$: the intersection over $f\in J$ of the projective tropical hypersurfaces
$V^{\trop}(f)$---that is, the $\semifield$-points of $\PP^n$ where the maximum in $f$ is attained at least twice---or the
set of $\semifield$-points of the closed subscheme defined by the bend congruence associated to $J$.  The definition of
the latter is introduced in \cite{GG1} and the equivalence between these two constructions is \cite[Proposition
5.1.6]{GG1}.

It is automatic that the tropical variety associated to a principal ideal $\langle f\rangle \subset \semifieldpolys$
equals the tropical hypersurface $V^{\trop}(f)$, since for any $fg\in \langle f\rangle$ we have
$V^{\trop}(f) \subset V^{\trop}(fg)$.  The following important result says that the same is also true for the Macaulay
tropical ideal $[f]$, at least when working over the tropical numbers.

\begin{theorem}
If $\semifield = \T$, then for any homogeneous $f\in \semifieldpolys$ the tropical variety $V^{\trop}([f])$ equals the tropical hypersurface $V^{\trop}(f)$.
\end{theorem}

\begin{remark}
We expect that this result holds for an arbitrary totally ordered idempotent semifield $\semifield$, not just the tropical semifield $\T$.
\end{remark}

\begin{proof}
To forestall confusion with our choice of $\max$ for the addition operation in~$\semifield$, 
we work in the max convention for $\T=\mathbb R\cup\{-\infty\}$ in this proof.

Since $f \in [f]$, what is to show is that $V^{\trop}(f) \subset V^{\trop}(h)$ for all $h\in [f]$.  By \cite[Theorem~6.3]{Fink-Rincon},
which assumes $\semifield=\T$,
for any fixed $h\in [f]$ there is a homogeneous $g\in \semifieldpolys$ 
such that $h$ is obtained from $gf$ by possibly decreasing certain
coefficients. The coefficients of $gf$ which may be decreased are
those in which the maximum in the computation of the product is obtained twice.  

Given a homogeneous element of $\semifieldpolys$,
say $k=\sum k_ax^a$ in multi-index notation, 
and a point $y\in \semifield^{n+1}$, let $P_y(k)$
be the convex hull within $\R^{n+1}$
of the indices $a$ such that $k_ay^a$ attains its maximum value.
Now fix $y \in V^{\trop}(f)$.
Then $P_y(f)$ is the convex hull of at least two points,
i.e., has at least two vertices.  
The same is true of $P_y(gf)$, which equals the Minkowski sum $P_y(g)+P_y(f)$.

The vertices of $P_y(gf)$ arise in a unique way as a point of $P_y(g)$
plus a point of $P_y(f)$.  Thus if $a$ is a vertex of $P_y(gf)$,
there is a unique $b$ such that $f_by^bg_{a-b}y^{a-b}$ takes its maximum value.
This $b$ is also then the unique index so that $f_bg_{a-b}$ takes its
maximum value.  This implies that the coefficient of $a$ in $h$ must equal
$f_bg_{a-b}$, not a lesser value; and it follows that every vertex
of $P_y(gf)$ is contained in $P_y(h)$.  In particular $P_y(h)$
contains at least two lattice points, so that $y\in V^{\trop}(h)$.
\end{proof}

\subsection{Weak order and the universal property of $[f]$}\label{sec:weak order}

We begin this section by showing that, at the level of underlying matroids,
$[f]$ is freest among the numerically principal ideals over~$f$.
Recall that given two matroids $M,M'$ on the same ground set~$E$, we say
$M'$ is a \emph{weak image} of~$M$ if every independent set of~$M'$ is independent in~$M$.
If $M$ and~$M'$ have the same rank, it is enough to check bases.

\begin{proposition}\label{prop:universal in weak order}
Let $f\in\semifieldpolys$ and let $I\subset\semifieldpolys$ be a numerically principal tropical ideal of degree $\deg f$ with $I_{\deg f}=\semifield f$.
For any degree $d$, the underlying matroid of~$I_d$ is a weak image of the underlying matroid of~$[f]_d$.
\end{proposition}

\begin{proof}
Suppose $K\subset\Delta^n_d$ with $|K|=N^n_{d-\deg f}$ is not a basis of~$[f]_d$ (by which we mean of its underlying matroid).
We must show that $K$ is also not a basis of~$I_d$.
In the proof we use $|$ and $/$ to denote the restriction and contraction operations on valuated matroids.

By Proposition~\ref{prop:transversal}, the underlying matroid of $[f]_d$ is the transversal matroid presented by the support of~$\D_d(f)$. 
Since $K$ admits no transversal, Hall's marriage theorem provides subsets $A\subseteq\Delta^n_{d-\deg f}$ and $B\subseteq K$ such that $|A|+|B|>N^n_{d-\deg f}$ and 
$\D_d(f)_{XY}=0$ for all $X\in A$ and $Y\in B$.

Because $I$ is a tropical ideal, $I_d$ contains $mf$ for all $m\in\Delta^n_{d-\deg f}$.
The multiples $mf$ for $m\in A$ lie in the intersection of $I_d$ with the coordinate subspace $I_d\cap\semifield^{\Delta^n_d\setminus B}$, which is a tropical linear space.
By Lemma~\ref{lem:monomial multiples rank bound},
\[\rank(I_d/B)=\rank(I_d\cap\semifield^{\Delta^n_d\setminus B})\ge|A|.\]
Because $I$ is numerically principal, $\rank I_d=N^n_{d-\deg f}$, so
\[
\rank(I_d|B)
= \rank(I_d)-\rank(I_d/B)
\le N^n_{d-\deg f}-|A|
< |B|.
\]
That is, $B$ is a dependent set in~$I_d$. So its superset $K$ is not a basis, as desired.
\end{proof}

\subsection{Comparison with the Maclagan-Rinc\'on ideal}

Maclagan and Rinc\'on observed that a family of matroids studied in \cite{Ardila-Billey} assemble to form a non-realizable homogeneous tropical ideal $I_{\mathrm{MR}} \subset \Bpolys$, for any $n \ge 2$ \cite[Example 2.8]{Maclagan-Rincon2}.  The Hilbert function of this ideal is that of a line, and when passing to $\T$ the associated tropical variety $V^{\trop}(I_{\mathrm{MR}}) \subset \PP^n(\T)$ is the standard tropical line.

In this section we review the construction of~$I_{\mathrm{MR}}$ and show that when $n=2$ it coincides with Macaulay ideal $[x+y+z]$.
The reader comparing to \cite{Maclagan-Rincon2} should bear Remark~\ref{rem:conventions} in mind: the matroids we associate in this paper to~$I_{\mathrm{MR}}$ are the dual of theirs.
In this section, for ease of translating to their conventions, 
we speak of our matroids through their cobases and cocircuits.


Let $M_d$ denote the matroid on $\Delta^n_d$ whose cobases are the $(d+1)$-element subsets $X$ of the simplex $\Delta^n_d$ that satisfy the following density bound:  for any $k\in \{0, \ldots, d\}$, each subsimplex $K \in \Delta^n_{k\to d}$ contains at most $k+1$ elements of $X$.  (One should think of this $k+1$ as the Hilbert function of a line.)

The cocircuits of $M_d$, which span $(I_{\mathrm{MR}})_d$ as a $\mathbb{B}$-module, are the subsets $C \subset \Delta^n_d$ that are contained in a
subsimplex $K\in \Delta^n_{(|C|-2) \to d}$  and such that for any $k < |C|-2$ each subsimplex $L \in \Delta^n_{k\to d}$ contains at most $k+1$
elements of $C$.  This says that the density bound for cobases is violated on $K$ and not on any smaller subsimplex.  From this description one immediately sees that if $C$ is a cocircuit in degree $d$ then $x_i C$ is a cocircuit in degree $d+1$, and so this sequence of matroids indeed constitutes a tropical ideal.

\begin{example}
Suppose $n=2$.  In degree 1, the cobases of $M_1$ are the following sets, 
drawn with black circles for monomials that are present:
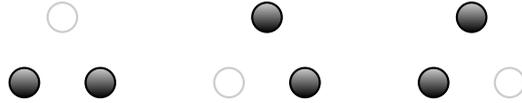
\begin{center}
\vspace{0.5\baselineskip}
\begin{tikzpicture}[font=\small]
 \path 
 ( 0,0) node [bdot] {}
 ( 1,0) node [bdot] {}
 ( 0.5,0.866) node [wdot] {};
\end{tikzpicture}
\hspace{10mm}
\begin{tikzpicture}[font=\small]
 \path 
 ( 0,0) node [wdot] {}
 ( 1,0) node [bdot] {}
 ( 0.5,0.866) node [bdot] {};
\end{tikzpicture}
\hspace{10mm}
\begin{tikzpicture}[font=\small]
 \path 
 ( 0,0) node [bdot] {}
 ( 1,0) node [wdot] {}
 ( 0.5,0.866) node [bdot] {};
\end{tikzpicture}
\end{center}
There is a single cocircuit consisting of all three monomials.

In degree 2, the cobases are the following sets plus all of their reflections and rotations:
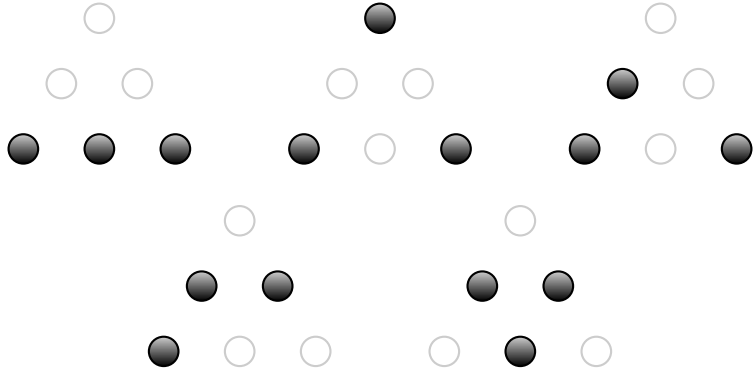
\begin{center}
\vspace{0.5\baselineskip}
\begin{tikzpicture}[font=\small]
 \path 
 ( 0,0) node [bdot] {}
 ( 1,0) node [bdot] {}
 ( 2,0) node [bdot] {}
 ( 0.5,0.866) node [wdot] {}
 ( 1.5,0.866) node [wdot] {}
 ( 1,1.732) node [wdot] {};
\end{tikzpicture}
\hspace{10mm}
\begin{tikzpicture}[font=\small]
 \path 
 ( 0,0) node [bdot] {}
 ( 1,0) node [wdot] {}
 ( 2,0) node [bdot] {}
 ( 0.5,0.866) node [wdot] {}
 ( 1.5,0.866) node [wdot] {}
 ( 1,1.732) node [bdot] {};
\end{tikzpicture}
\hspace{10mm}
\begin{tikzpicture}[font=\small]
 \path 
 ( 0,0) node [bdot] {}
 ( 1,0) node [wdot] {}
 ( 2,0) node [bdot] {}
 ( 0.5,0.866) node [bdot] {}
 ( 1.5,0.866) node [wdot] {}
 ( 1,1.732) node [wdot] {};
\end{tikzpicture}

\vspace{0.5cm}
\begin{tikzpicture}[font=\small]
 \path 
 ( 0,0) node [bdot] {}
 ( 1,0) node [wdot] {}
 ( 2,0) node [wdot] {}
 ( 0.5,0.866) node [bdot] {}
 ( 1.5,0.866) node [bdot] {}
 ( 1,1.732) node [wdot] {};
\end{tikzpicture}
\hspace{10mm}
\begin{tikzpicture}[font=\small]
 \path 
 ( 0,0) node [wdot] {}
 ( 1,0) node [bdot] {}
 ( 2,0) node [wdot] {}
 ( 0.5,0.866) node [bdot] {}
 ( 1.5,0.866) node [bdot] {}
 ( 1,1.732) node [wdot] {};
\end{tikzpicture}
\end{center}
The cocircuits are given by all reflections and rotations of the three following sets. 
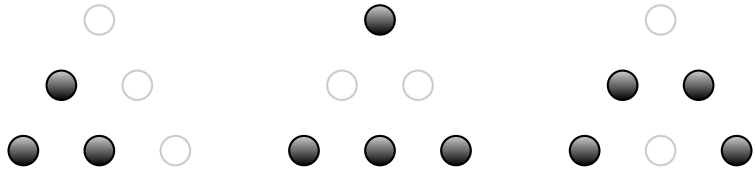
\begin{center}
\vspace{0.5\baselineskip}
\begin{tikzpicture}[font=\small]
 \path 
 ( 0,0) node [bdot] {}
 ( 1,0) node [bdot] {}
 ( 2,0) node [wdot] {}
 ( 0.5,0.866) node [bdot] {}
 ( 1.5,0.866) node [wdot] {}
 ( 1,1.732) node [wdot] {};
\end{tikzpicture}
\hspace{10mm}
\begin{tikzpicture}[font=\small]
 \path 
 ( 0,0) node [bdot] {}
 ( 1,0) node [bdot] {}
 ( 2,0) node [bdot] {}
 ( 0.5,0.866) node [wdot] {}
 ( 1.5,0.866) node [wdot] {}
 ( 1,1.732) node [bdot] {};
\end{tikzpicture}
\hspace{10mm}
\begin{tikzpicture}[font=\small]
 \path 
 ( 0,0) node [bdot] {}
 ( 1,0) node [wdot] {}
 ( 2,0) node [bdot] {}
 ( 0.5,0.866) node [bdot] {}
 ( 1.5,0.866) node [bdot] {}
 ( 1,1.732) node [wdot] {};
\end{tikzpicture}
\end{center}
The cocircuits in the orbit of the left configuration violate the cobasis density bound in a simplex $K \in \Delta^2_{1\to 2}$, and the other two orbits violate the density bound only on the whole of $\Delta^2_2$, as they consist of all sets of four monomials that don't contain one of the cocircuits of the left type.
\end{example}

\begin{proposition}\label{prop:MR}
When $n=2$, the Maclagan-Rinc\'on tropical ideal $I_{\mathrm{MR}} \subset S=\mathbb{B}[x_0, x_1, x_2]$ is equal to $[x_0
+ x_1 + x_2]$.
\end{proposition}

\begin{proof}
Theorem 6.3 of \cite{Ardila-Billey} asserts that the matroids $M_d$ of
$I_{\mathrm{MR}}$ are transversal, and the first proof they give
shows that there is a presentation identical to the canonical
presentation of $[x_0 + x_1 + x_2]_d$.
\end{proof}

In the language of the proof of Proposition~\ref{prop:universal in weak order},
if a set $X$ fails to be a cobasis of~$M_d$ 
because it violates the density bound on a subsimplex $m\Delta^2_\ell \in \Delta^2_{\ell\to d}$,  then the sets $A=m\Delta^2_{\ell-1}$ and $B=(\Delta^2_d\setminus m\Delta^2_\ell)\setminus X$
provide by Hall's marriage theorem an obstruction to the existence of a transversal to~$\Delta^2_d\setminus X$.
The content of \cite[Theorem 6.3]{Ardila-Billey} is that, in fact,
for any non-(co)basis of~$M_d$, the obstruction in Hall's marriage theorem may be chosen to be of this form.

It is tempting to try to generalize the construction of $I_{\mathrm{MR}}$ to give hypersurfaces of higher degree as follows. 
\begin{recipe}\label{recipe}
Propose $X\subset \Delta^n_d$ of size $N^n_d-N^n_{d-D}$ to be a cobasis if it satisfies a density bound given by the Hilbert function of a degree $D$ hypersurface: i.e., if each $K\in \Delta^n_{k\to d}$ satisfies $|X\cap K| \leq N^n_k - N^n_{k-D}$. 
\end{recipe}
In degree $D$ (the bottom degree), this recipe does indeed define the uniform matroid of rank 1 on $\Delta^n_D$, so the tropical linear space is the line spanned by $f=\sum_{m\in\Delta^2_D} m$ (the sum of all degree $D$ monomials).

For $D=2$, this might work.  Unfortunately, for $D\geq 3$, this recipe fails to yield a matroid in higher degrees.  This is because we do not have control over the Hall obstructions as we did for $D=1$, and therefore there is no generalization of Proposition~\ref{prop:MR} describing the cociruits of~$[f]$ only in terms of cardinalities of intersections with subsimplices.
The proof of the next proposition gives an explicit non-cobasis of $[f]$ for $D=3$ that cannot be detected using just these cardinalities. The example generalizes to $D>3$.

\begin{proposition}\label{prop:not just simplices}
For $n=2$ and $D=3$, the cobases proposed by Recipe \ref{recipe} fail to define a matroid in degree 6.
\end{proposition}
\begin{proof}
The proof is by contradiction.  Assume that Recipe \ref{recipe} does yield the cobases of a matroid in degree 6.  Consider $f=\sum_{m\in\Delta^2_3}m=x^3+x^2y+\cdots+z^3$ in $\mathbb{B}[x,y,z]$.  In degree $6$, the supports of $x^3f$ and $y^3f$ are both cocircuits (observe that removing any element yields a cobasis).  Applying the cocircuit elimination axiom to the unique common monomial $x^3y^3$ of these two sets shows that the support of
\[
g=x^3\left(\sum_{m\in\Delta^2_3\setminus\{y^3\}}m\right)
+ y^3\left(\sum_{m\in\Delta^2_3\setminus\{x^3\}}m\right),
\]
illustrated in Figure \ref{fig:not just simplices}, must be a union of cocircuits.   In particular,  $\supp(g)$ is not a cobasis.

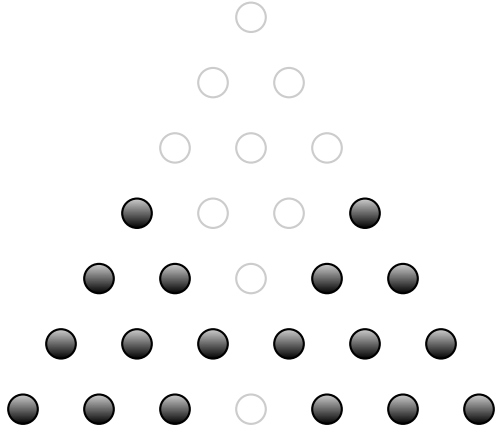
\begin{figure}[t]
    \centering
\begin{tikzpicture}[font=\small]
 \path 
 ( 0,0) node [bdot] {}
 ( 1,0) node [bdot] {}
 ( 2,0) node [bdot] {}
 ( 3,0) node [wdot] {}
 ( 4,0) node [bdot] {}
 ( 5,0) node [bdot] {}
 ( 6,0) node [bdot] {}
 ( 0.5,0.866) node [bdot] {}
 ( 1.5,0.866) node [bdot] {}
 ( 2.5,0.866) node [bdot] {}
 ( 3.5,0.866) node [bdot] {}
 ( 4.5,0.866) node [bdot] {}
 ( 5.5,0.866) node [bdot] {}
 ( 1,1.732) node [bdot] {}
 ( 2,1.732) node [bdot] {}
 ( 3,1.732) node [wdot] {}
 ( 4,1.732) node [bdot] {}
 ( 5,1.732) node [bdot] {}
 ( 1.5,2.598) node [bdot] {}
 ( 2.5,2.598) node [wdot] {}
 ( 3.5,2.598) node [wdot] {}
 ( 4.5,2.598) node [bdot] {}
 ( 2,3.464) node [wdot] {}
 ( 3,3.464) node [wdot] {}
 ( 4,3.464) node [wdot] {}
 ( 2.5,4.330) node [wdot] {}
 ( 3.5,4.330) node [wdot] {}
 ( 3,5.196) node [wdot] {};
\end{tikzpicture}
    \caption{Configuration from the proof of Proposition \ref{prop:not just simplices} showing that Recipe \ref{recipe} does not define a matroid in degree 6 for $n=2$ and $D=3$.  For $f=\sum_{m\in\Delta^2_3}m=x^3+x^2y+\cdots+z^3$, this set a not cobasis of the Macaulay tropical ideal $[f]$ or any tropical ideal containing $f$.  The monomials are arranged as in Figure~\ref{fig: Delta_d^2}.}
    \label{fig:not just simplices}
\end{figure}

On the other hand, let us check that $\supp(g)$ satisfies the density bound of Recipe \ref{recipe}, which says in this case that cobases are sets having at most $N^2_k - N^2_{k-3} = 3k$ elements lying in each subsimplex $K\in \Delta^n_{k\to d}$ for $k=1\ldots d$.  By inspection, one sees that this density bound is indeed satisfied by $\supp(g)$, and so Recipe \ref{recipe} says it is a cobasis.  We thus have the desired contradiction.
\end{proof}

\appendix
\section{Minimal circuit eliminations (by Joshua Mundinger)}

The (co)vector elimination characterization of tropical linear spaces mentioned in \S\ref{sec:TLS} suggests a naive method for expanding any subset of $\semifield^n$ to a tropical linear space, and hence of expanding any ideal in $\semifieldpolys$ to a tropical ideal, that we briefly explore here.

In what follows, for $v = \sum_i v_i e_i \in \semifield^n$ we write  
\[v_{\widehat{j}} = \sum_{i \ne j} v_ie_i.\]  

\begin{adefinition}
We say that $w = \sum_i w_ie_i \in \semifield^n$ is a $j$-\emph{elimination} of a pair \[v = \sum_i v_ie_i,~v' = \sum_i v'_ie_i\] in $\semifield^n$ if:
\begin{enumerate}
\item $v_j = v'_j \ne 0$ and $w_j = 0$; 
\item $w+v+v' = v+v'$; and, 
\item $w_i = v_i + v'_i$ for all indices $i$ such that $v_i \ne v'_i$. 
\end{enumerate}
We call $w$ a \emph{minimal} $j$-elimination of the pair $v,v'$ if $v_j=v'_j \ne 0$ and \[w = (v+v')_{\widehat{j}}\] (this implies the above three conditions).
\end{adefinition}

\begin{aremark}
It follows from the (co)vector characterizations of valuated matroids that a submodule $M \subset \semifield^n$ is a tropical linear space if and only if the following holds: whenever $v,v'\in M$ satisfy $v_j = v'_j \ne 0$ for some index $j$, then $M$ contains a $j$-elimination of the pair $v,v'$. 
\end{aremark}

\begin{adefinition}
For any subset $M \subset \semifield^n$, let $M^{[0]}$ be the submodule generated by $M$ and inductively define $M^{[\ell]} \subset \semifield^n$ to be the submodule generated by $M^{[\ell-1]}$ and by the collection of  minimal eliminations $(v+v')_{\widehat{j}}$ of pairs $v,v'\in M^{[\ell-1]}$ such that (i) $v_j = v'_j \ne 0$ and (ii) $M^{[\ell-1]}$ does not contain a $j$-elimination of the pair $v,v'$.  We denote the infinite union by \[M^{[\infty]} = \bigcup_{\ell = 0}^\infty M^{[\ell]}\]
and call it the \emph{minimal-elimination envelope} of $M$.
\end{adefinition}

Since these modules are nested,
\[M^{[0]} \subset M^{[1]} \subset M^{[2]} \subset \cdots,\]
the minimal-elimination envelope $M^{[\infty]} \subset \semifield^n$ is a submodule containing $M$. 

\begin{aproposition}\label{prop:minenvTLS}
The submodule $M^{[\infty]} \subset \semifield^n$ is a tropical linear space.
\end{aproposition}

\begin{proof}
For any pair $v,v'\in M^{[\infty]}$ and any index $j$ such that $v_j = v'_j \ne 0$, we must have $v\in M^{[\ell]}$ and $v'\in M^{[\ell']}$ for some integers $\ell,\ell'$.  Without loss of generality, assume that $\ell \le \ell'$.  Then $v \in M^{[\ell']}$ since $M^{[\ell]} \subset M^{[\ell']}$, so either $M^{[\ell']}$ contains a $j$-elimination of the pair $v,v'$ or else by definition $M^{[\ell'+1]}$ contains a minimal $j$-elimination of this pair.  In either case, $M^{[\infty]}$ contains a $j$-elimination of the pair.
\end{proof}

We now use this minimal-elimination envelope construction degree-by-degree to create a procedure for extending a homogeneous ideal to a tropical ideal.

\begin{adefinition}
Let $I \subset \semifieldpolys$ be a homogeneous ideal.  We define the \emph{minimal-elimination envelope} $[I]^{\infty} \subset \semifieldpolys$ as the direct sum of the following $\semifield$-modules: $[I]_0^{\infty} = I_0,~[I]_1^{\infty} = I_1^{[\infty]}$, and $[I]_d^{\infty}$ for $d > 1$ is $M^{[\infty]}$ where $M \subset \semifieldpolys_d$ is the $\semifield$-module generated by $I_d$ and by the submodules $x_i[I]_{d-1}^{\infty}$ for each $i=0,\ldots,n$.
\end{adefinition}

This $\semifield$-module $[I]^{\infty}$ is a homogeneous ideal since it is graded and satisfies $x_i[I]^{\infty} \subset [I]^{\infty}$ for each $i=0,\ldots,n$.  Moreover, it is a tropical ideal since each graded piece is a tropical linear space by Proposition \ref{prop:minenvTLS}.  The problem with this construction is that we have no control over the Hilbert function, and it may even result in the irrelevant ideal.  However, in one very special case it behaves nicely and in fact coincides with the Macaulay tropical ideal construction introduced in this paper:

\begin{atheorem} \label{theorem: main}
	Let $f = a_\alpha x^\alpha + a_\beta x^\beta \in \semifieldpolys$ be a homogeneous binomial.  
 Then the ideal $[\langle f \rangle]^{\infty}$ is exactly the Macaulay tropical ideal $[f]$, which is also equal to 
 $\langle f\rangle^{\perp \perp}$ in this case.
\end{atheorem}

\begin{aremark}
In general the minimal elimination envelope for a principal ideal does \emph{not} coincide with the Macaulay tropical ideal.  Indeed, symbolic calculations performed in SageMath show that for the linear trinomial $ f = x_0 + x_1+ x_2 \in \Bpolys$, the tropical linear spaces $[\langle f \rangle]^{\infty}_3$ and $[f]_3$ are unequal, although curiously they do have the same rank.
\end{aremark}

Before proving the preceding theorem, we establish a few lemmas.

\begin{alemma} \label{lemma: transversal-trees}
	Let $M$ be a transversal matroid on $\{1,\ldots,n\}$ 
	presented by $n-1$ sets $A_1, \ldots, A_{n-1}$, all of size two.
	Let $G$ be the graph with vertices $\{1,\ldots,n\}$ and edges $\{A_1, \ldots, A_{n-1}\}$.
	If $G$ is a tree, then $M = U_{n, n-1}$ is the rank $n-1$ uniform matroid.
\end{alemma}
\begin{proof}
	The proof is by induction on $n$.
	We need to show that for any $1 \le i \le n$, the set
	$\{1,\ldots,n\}-i$ is a transversal of $\{A_1, \ldots, A_{n-1}\}$.
	Since $G$ is connected, there is an edge containing $i$, say $A_{n-1} = \{i,j\}$.
	Let $G' = G - A_{n-1}$, which is a disjoint union of two trees since $G$ is a tree.
	One connected component of $G'$ contains $i$, and the other contains $j$;
	let $V_i$ denote the vertices in the component containing $i$ and similarly for $V_j$.
	By inductive hypothesis, $V_i - i$ is a transversal of the sets $A_k$ contained in $V_i$,
	and $V_j - j$ is a transversal of the sets $A_k$ contained in $V_j$;
	thus, $\{1,\ldots,n\} - \{i,j\}$ is a transversal of $A_1, \ldots, A_{n-2}$.
	Since $A_{n-1} = \{i,j\}$, we conclude that $\{1,\ldots,n\} - i$ is a transversal of $A_1, \ldots, A_{n-2}, A_{n-1}$.
\end{proof}

\begin{alemma}\label{lem:Aequivrel}
	Let $L\subseteq \B^n$ be a tropical linear space.
	The relation on $\{1,\ldots,n\}$ defined by $i \sim j$ if $i = j$ or $e_i + e_j \in L$
	is an equivalence relation.
\end{alemma}
\begin{proof}
	Reflexivity and symmetry are by definition, transitivity is by vector elimination:
	if $i,j,k$ are all distinct and $e_i + e_j, e_j + e_k \in L$,
	then $e_i + e_k \in L$ as well.
\end{proof}

\begin{alemma} \label{lemma: elimination}
	Let $A \subseteq \B^n$ be a submodule generated by elements with support of size two.
	Let $\sim$ be the equivalence relation generated by $i \sim j$ if $e_i + e_j \in A$.
	Then $A^{[\infty]} = \langle e_i + e_j: i \sim j, i \neq j \rangle$.
\end{alemma}
\begin{proof}
	Let $B = \langle e_i + e_j~|~i \sim j,~i\neq j\rangle$.  Since $A^{[\infty]}$ is a tropical linear space, it must contain $\{e_i + e_j~|~ i \sim j,~i \neq j\}$ by Lemma \ref{lem:Aequivrel}, and so $B \subset A^{[\infty]}$.  Thus we must establish the reverse containment.  Clearly $A \subset B$, so it suffices to show that $B$ is closed under minimal elimination.  Let $P$ be the partition of $\{1,\ldots,n\}$ of equivalence classes under $\sim$; then $B$ is the tropical linear space associated to the matroid direct sum of uniform matroids $\oplus_{E \in P} U_{E, |E|-1}$.   The property of being closed under minimal elimination is closed under matroid direct sum, so it suffices to show that the tropical linear space associated to any uniform matroid $U_{n,d}$ is closed under minimal elimination; but this is clear since $L_{U_{n,d}}$ comprises all vectors with support of size at least $n-d+1$ and $|\supp (v+v')_{\widehat{j}}| \geq |\supp(v)|$.
\end{proof}

\begin{proof}[Proof of Theorem \ref{theorem: main}]
We begin by reducing to the case $f = x^\alpha + x^\beta$ with $\semifield=\B$. 
We may embed $\semifield^\times$ in a divisible ordered abelian group $\Gamma$, and thus $\semifield$ in $\semifield'=(\Gamma\cup\{0\},\max,\cdot)$.
Divisibility implies that we can rescale monomials
so that the binomial $f$ becomes $x^\alpha + x^\beta$, defined over $\B\subset\semifield'$. 
Both constructions $[\langle x^\alpha + x^\beta\rangle]^\infty$ and $[x^\alpha + x^\beta]$ 
are then the tensor product with $\semifield'$ of the constructions for~$\B$,
and both are defined over $\semifield$ once the rescaling is undone.

	Suppose $f = x^\alpha + x^\beta$ is of degree $d$, 
	and let $e \ge d$.
	Let $\sim_e$ be the equivalence relation $\Delta_e^n$  
	generated by $x^{\alpha + \gamma} \sim_e x^{\beta + \gamma}$ for all $x^\gamma \in \Delta_{e-d}^n$.
	We claim that \[[\langle f \rangle]_e^\infty = \langle x^u + x^v~|~x^u \sim_e x^v,~u \neq v\rangle = [f]_e.\]

The proof for $[\langle f \rangle]_e^\infty$ is by induction on $e$.  The base case $e=d$ is clear, so assume $e > d$.  Since \[\langle f \rangle_e = \sum_{i=0}^n x_i\langle f \rangle_{e-1} \subset \sum_{i=0}^n x_i[\langle f\rangle]^\infty_{e-1},\] the module $[\langle f \rangle]_e^\infty$ is the minimal-elimination envelope of the modules $x_i[\langle f\rangle]^\infty_{e-1}$ for all $i=0,\ldots,n$.  By the inductive hypothesis and Lemma \ref{lemma: elimination}, $[\langle f \rangle]_{e-1}^\infty$ is the minimal-elimination envelope of the binomials $x^{\alpha+\lambda} + x^{\beta+\lambda}$ for $x^\lambda\in\Delta^n_{e-1-d}$, so $x_i[\langle f \rangle]_{e-1}^\infty$ is the minimal-elimination envelope of the binomials $x_ix^{\alpha+\lambda} + x_ix^{\beta+\lambda}$ for $x^\lambda\in\Delta^n_{e-1-d}$, hence $[\langle f \rangle]_e^\infty$ is the minimal-elimination envelope of the binomials $x^{\alpha+\gamma} + x^{\beta+\gamma}$ for $x^\gamma\in\Delta^n_{e-d}$.  The result then follows from Lemma \ref{lemma: elimination}.

	Now for $[f]_e$, the tropical linear space associated to the transversal matroid with presentation \[\{\{x^{\alpha +\gamma},~x^{\beta + \gamma}\}~|~ x^\gamma \in \Delta^n_{e-d}\}.\]
	Let $G$ be the graph on vertices $\Delta^n_e$ with edges given by the above presentation; the connected components of $G$ are the equivalence classes of our equivalence relation $\sim_e$.  We claim that $G$ is a forest.  Since $x^\alpha \neq x^\beta$,
	there is an index $i$ such that $\alpha_i < \beta_i$.
	Since each edge of $G$ is of the form $\{x^{\alpha +\gamma},x^{\beta + \gamma}\}$,
	any monomial $x^\delta$ is incident to at most the vertices $x^{\delta - \alpha + \beta}$
	and $x^{\delta + \alpha - \beta}$.
	But \[(\delta - \alpha + \beta)_i > \delta_i > (\delta + \alpha - \beta)_i,\]
	so any path has a strictly monotonic power of $x_i$ and hence $G$ has no cycles.  We can now apply Lemma \ref{lemma: transversal-trees} to the connected components of $G$, which shows that $[f]_e$ is the tropical linear space associated to the direct sum, over the equivalence classes $E$ of   $\sim_e$, of uniform matroids $U_{E,|E|-1}$.  This is the module generated by $x^u+x^v$ for $x^u\sim_e x^v,~u\ne v$, as claimed.
\end{proof}

\bibliographystyle{amsalpha}
\bibliography{bib}

\providecommand{\bysame}{\leavevmode\hbox to3em{\hrulefill}\thinspace}
\providecommand{\MR}{\relax\ifhmode\unskip\space\fi MR }
\providecommand{\MRhref}[2]{%
  \href{http://www.ams.org/mathscinet-getitem?mr=#1}{#2}
}
\providecommand{\href}[2]{#2}
\begin{thebibliography}{Mun18}

\bibitem[AB07]{Ardila-Billey}
Federico Ardila and Sara Billey, \emph{Flag arrangements and triangulations of
  products of simplices}, Adv. Math. \textbf{214} (2007), no.~2, 495--524.

\bibitem[AR22]{Anderson-Rincon-paving}
Nicholas Anderson and Felipe Rinc{\'o}n, \emph{Paving tropical ideals}, J.
  Algebr. Comb. \textbf{56} (2022), no.~1, 101--116 (English).

\bibitem[Bon10]{BoninTransversalNotes}
Joseph~E Bonin, \emph{An introduction to transversal matroids},
  \url{https://maa.org/sites/default/files/pdf/shortcourse/2011/TransversalNotes.pdf},
  2010.

\bibitem[DW91]{Dress-Wenzel-1}
Andreas Dress and Walter Wenzel, \emph{Grassmann-{P}l\"ucker relations and
  matroids with coefficients}, Adv. Math. \textbf{86} (1991), no.~1, 68--110.

\bibitem[DW92]{Dress-Wenzel-2}
Andreas W.~M. Dress and Walter Wenzel, \emph{Valuated matroids}, Adv. Math.
  \textbf{93} (1992), no.~2, 214--250.

\bibitem[EF65]{Edmonds-Fulkerson}
Jack Edmonds and D.~R. Fulkerson, \emph{Transversals and matroid partition}, J.
  Res. Nat. Bur. Standards Sect. B \textbf{69B} (1965), 147--153.

\bibitem[FR15]{Fink-Rincon}
Alex Fink and Felipe Rinc{\'o}n, \emph{Stiefel tropical linear spaces}, J.
  Combin. Theory Ser. A \textbf{135} (2015), 291--331.

\bibitem[FR16]{Dhruv-Tyler-Hahn}
Tyler Foster and Dhruv Ranganathan, \emph{Hahn analytification and connectivity
  of higher rank tropical varieties}, Manuscripta Math. \textbf{151} (2016),
  no.~3-4, 353--374. \MR{3556824}

\bibitem[Fre13]{Frenk}
Bart Frenk, \emph{Tropical varieties, maps and gossip}, Ph.D. thesis, Eindhoven
  University of Technology, 2013, \url{https://doi.org/10.6100/IR750815}.

\bibitem[GG16]{GG1}
Jeffrey Giansiracusa and Noah Giansiracusa, \emph{Equations of tropical
  varieties}, Duke Math. J. \textbf{165} (2016), no.~18, 3379--3433.

\bibitem[GG18]{GG3}
\bysame, \emph{A {G}rassmann algebra for matroids}, Manuscripta Math.
  \textbf{156} (2018), no.~1-2, 187--213.

\bibitem[Lor22]{Lorscheid-hyperfields}
Oliver Lorscheid, \emph{Tropical geometry over the tropical hyperfield}, Rocky
  Mt. J. Math. \textbf{52} (2022), no.~1, 189--222 (English).

\bibitem[Lor23]{Lorscheid-unifying}
\bysame, \emph{A unifying approach to tropicalization}, Trans. Am. Math. Soc.
  \textbf{376} (2023), no.~5, 3111--3189 (English).

\bibitem[MR18]{Maclagan-Rincon2}
Diane Maclagan and Felipe Rinc\'on, \emph{Tropical ideals}, Compos. Math.
  \textbf{154} (2018), no.~3, 640--670.

\bibitem[MR20]{Maclagan-Rincon}
\bysame, \emph{Tropical schemes, tropical cycles, and valuated matroids}, J.
  Eur. Math. Soc. (JEMS) \textbf{22} (2020), no.~3, 777--796 (English).

\bibitem[MR22]{Maclagan-Rincon3}
\bysame, \emph{Varieties of tropical ideals are balanced}, Adv. Math.
  \textbf{410} (2022), 44 (English), Id/No 108713.

\bibitem[MS15]{Maclagan-Sturmfels}
Diane Maclagan and Bernd Sturmfels, \emph{Introduction to tropical geometry},
  Graduate Studies in Mathematics, vol. 161, American Mathematical Society,
  Providence, RI, 2015.

\bibitem[MT01]{Murota-circuits}
Kazuo Murota and Akihisa Tamura, \emph{On circuit valuation of matroids}, Adv.
  in Appl. Math. \textbf{26} (2001), no.~3, 192--225.

\bibitem[Mun18]{Mundinger}
J.~Mundinger, \emph{The image of a tropical linear space}, arXiv:1808.02150,
  2018.

\bibitem[PW70]{Piff-Welsh}
M.~J. Piff and D.~J.~A. Welsh, \emph{On the vector representation of matroids},
  J. London Math. Soc. (2) \textbf{2} (1970), 284--288.

\bibitem[Sil22]{Silversmith}
R.~Silversmith, \emph{The matroid stratification of the {Hilbert} scheme of
  points on {{\(\mathbb{P}^1\)}}}, Manuscr. Math. \textbf{167} (2022), no.~1-2,
  173--195 (English).

\bibitem[Spe08]{Speyer}
David Speyer, \emph{Tropical linear spaces}, SIAM J. Discrete Math. \textbf{22}
  (2008), no.~4, 1527--1558.

\end{thebibliography}

\end{document}